\documentclass{amsart}
\copyrightinfo{2014}{the author}

\usepackage{amsmath}
\usepackage{amsfonts}
\usepackage{amsthm}
\usepackage{amssymb}
\usepackage{graphicx, url, todonotes}

\newtheorem{theorem}{Theorem}[section]

\newtheorem{lemma}[theorem]{Lemma}

\theoremstyle{definition}
\newtheorem{definition}[theorem]{Definition}

\newtheorem*{Maekawa}{Maekawa's Theorem}

\title{Coloring Connections with Counting Mountain-valley Assignments}

 \author[Hull]{Thomas C. Hull}
 \address{Hull: Department of Mathematics, Western New England University, Springfield MA, USA}
\email{thull@wne.edu}
\thanks{Appeared in {\em Origami$^6$}, Miura et. al. eds., The American Mathematical Society, pp. 3-10.}

\subjclass{Primary 05C15, 05C30; Secondary  52C99}

\keywords{flat folding, graph coloring, Miura-ori}

\begin{document}

\maketitle

%\begin{abstract}
%We survey more recent attempts at enumerating the number of mountain-valley assignments that allow a given crease pattern to locally fold flat.  In particular, we solve this problem for square twist tessellations and generalize the method used to a broader family of crease patterns.  We also describe the more difficult case of the Miura-ori and a recently-discovered bijection with 3-vertex colorings of grid graphs.
%\end{abstract}

\section{Introduction}

In the history of origami-mathematics, studying the combinatorics of flat folds has been consistently difficult.  The stamp-folding problem, where we aim to count the number of different ways to fold a grid of postage stamps, might be the oldest such problem which is still unsolved \cite{U}. Determining the number of ways a given crease pattern can fold flat has also been found to have applications in science, such as in enumerating the energy states of a crumpled polymer sheet or in determining the mechanics of self-folding materials.  

Many such problems amount to counting the number of valid mountain-valley (MV) assignments that will fold a given crease pattern into a flat state.  Here {\em valid} means that it will allow the crease pattern to fold flat without the paper self-intersecting or ripping.  A recent advancement in this area is to convert such MV-assignment enumeration problems to a graph coloring problem and then use the vast field of graph colorings to help solve the problem.  

Note that this general approach is not new;  physicists have been converting folding problems from polymer materials science to coloring problems for over a decade \cite{DiF}.  In such work, however, the folding problems are typically counting the number of foldings whose crease patterns are a subset of a larger grid of creases, not focusing on a specific crease pattern and counting MV-assignments.  In this paper we briefly survey recent approaches in linking the enumeration of valid MV-assignments to graph coloring problems.

\section{Preliminaries}

To begin we briefly review a few basic results from the combinatorics of flat origami.  See \cite{Hull2} for a more detailed treatment.

Consider our crease pattern $C$ to be a planar graph $C=(V,E)$ drawn on a piece of paper.  We partition the vertices into two classes, $V=V_B\cup V_I$ where $V_B$ are the vertices on the boundary of the paper and $V_I$ are the vertices in the paper's interior.  The edges in $E$ are the crease lines of our origami model, and in this paper we will only consider crease patterns $C$ that {\em fold flat}, or are {\em flat-foldable}, which means that all the creases can be folded to transform (in a piecewise isometric manner) the paper into a flat object such that the paper does not self-intersect.

We define a {\em MV-assignment} for a crease pattern $C=(V,E)$ to be a function $\mu:E\rightarrow\{-1,1\}$.  Here $-1$ indicates that the crease is a valley and $1$ indicates a mountain crease.  A MV-assignment $\mu$ is called {\em valid} if the crease pattern can fold flat (without self-intersecting) using the mountains and valleys as specified by $\mu$.

Many results on flat-foldable crease patterns are {\em local} in that they pertain only to the foldability or MV-assignment at a single interior vertex.  One of the most fundamental such results is Maekawa's Theorem.

\begin{Maekawa}
The difference between the number of mountain and valley creases meeting at a vertex in a flat-foldable crease pattern is two. That is, for the creases $c_i$ adjacent to the vertex we must have $\sum \mu(c_i) = \pm2$.
\end{Maekawa}

See \cite{Hull1}, \cite{Hull2}, or \cite{H1} for a proof.  This means at a degree 4 vertex in a flat-foldable crease pattern, we must have 3 mountain and 1 valley creases or vice-versa.  Note that Maekawa's Theorem is only a necessary condition for a vertex to fold flat.  

Sometimes specifying the angles between the creases meeting at a vertex is needed.  If the creases meeting at a vertex in a crease pattern are denoted $l_i$ for $i=1,...,n$, then we let $\alpha_i$ denote the angles, in order, between these creases, where $\alpha_i$ is between creases $l_i$ and $l_{i+1}$.  The angles $\alpha_i$ are called the {\em angle sequence} for the vertex.

Another fundamental consideration is the fact that local conditions in flat-foldability do not typically translate to global conditions.  The biggest such obstacle is the fact that if a crease pattern $C$ with a MV-assignment $\mu$ satisfies the condition that every vertex locally folds flat under $\mu$, then it is not guaranteed to globally fold flat.  Examples of crease patterns that are locally but not globally flat-foldable can be found in \cite{Hull1, Hull2, H1}.  The best attempt in print at trying to create conditions for global flat-foldability can be found in \cite{Justin}.  In 1996 Bern and Hayes proved that deciding whether or not a given crease pattern can globally fold flat is NP-hard, even if the MV-assignment is given as well \cite{Bern}.  

For this reason, when counting valid MV-assignments we usually only count those that locally fold flat.  Adding the global constraint makes such problems much more difficult, and there are currently no general techniques for addressing the global flat-foldability problem in MV-assignment enumeration.

\section{Two-colorable crease patterns}

The simplest cases of MV-assignments in flat-foldable crease patterns are those which are equivalent to properly 2-vertex coloring a graph.  These are crease patterns where the mountains and valleys can be determined using simple rules.  One of the most basic such rules in the world of flat-foldability is the following (see \cite{DOR07}):%(see \cite{Hull2, Hull3}):

\begin{lemma}[Big-Little-Big Lemma]\label{lem:blb}
Let $v$ be a flat vertex fold with angle sequence $\alpha_i$ and a valid MV assignment $\mu$.  If $\alpha_{i-1}>\alpha_i<\alpha_{i+1}$ for some $i$, then $\mu(l_i)\not= \mu(l_{i+1})$.
\end{lemma}

The proof of this Lemma is simple; if $\mu(l_i)=\mu(l_{i+1})$, then the sector of paper made by angle $\alpha_i$ will be covered on the same side of the paper by the two bigger adjacent angle sectors, forcing the paper to self-intersect or to create additional creases.

The Big-Little-Big Lemma allows us to examine crease patterns and quickly determine if some creases are forced to be the same or different.  It can't be applied to all flat-foldable crease patterns, but it can be a great aid in certain families of crease patterns.  In fact, we can make the influence of the Big-Little-Big Lemma and Maekawa's Theorem more precise by thinking of  a MV-assignment $\mu$ as a 2-coloring of the edges (with colors $-1$ and $+1$).   Our goal, then, is to convert our MV-assignment into a proper 2-vertex coloring of a graph.

Before we do this, recall that in graph theory $P_2$ denotes a {\em path of length 2}, which has three vertices $a$, $b$, and $c$ and edges $\{a,b\}$ and $\{b,c\}$.  The vertices $a$ and $c$ are called the {\em ends of $P_2$} and must get the same color when we 2-color the vertices of $P_2$.  In fact, the end vertices must always get the same color when 2-coloring the vertices of any even-length path $P_{2n}$.

\begin{definition}\label{def:origami line graph}
Given a flat-foldable crease pattern $C=(V,E)$, we define the {\em origami line graph} $C_L=(V_L,E_L)$ to be a graph produced as follows:  Start with an initial set of vertices $V_L$ to be the creases $\{c_1, ..., c_n\}$ in $E$.  Then perform the following steps:
\begin{enumerate}
\setlength{\itemsep}{0pt}
\item For all pair of creases $c_i, c_j\in E$, if they are forced to have different MV parity, then let $\{c_i,c_j\}\in E_L$.
\item For all pair of creases $c_i, c_j\in E$, if they are forced to have the same MV parity and $c_i$ and $c_j$ are not already the ends of a path of even length from performing step 1, then add a new vertex $v_{i,j}$ to $V_L$ and let $\{c_i, v_{i,j}\}, \{v_{i,j},c_j\}\in E_L$. 
\end{enumerate}
\end{definition}

This definition is an extension of one found in \cite{Hull1}.  The point of the origami line graph is the following:

\begin{theorem}\label{thm:lg}
Given a crease pattern $C$, if the origami line graph $C_L$ is not properly 2-vertex colorable, then $C$ is not flat-foldable.
\end{theorem}

For the purposes of counting, if a crease pattern has the property that all of the ways in which mountains and valleys can influence each other is captured by the origami line graph, then counting the number of valid MV-assignments is just a matter of counting the number of proper 2-colorings of the origami line graph.  We say that in a flat-foldable crease pattern $C$, the MV-assignments are {\em determined by $C_L$} if every MV-assignment $\mu$ corresponds to a unique 2-vertex coloring of $C_L$ and vice-versa.    There are only two ways to properly 2-color the vertices of a connected, 2-colorable graph, and therefore we have the following:

\begin{theorem}\label{thm:2color}
Let $C$ be a flat-foldable crease pattern whose MV-assignments are determined be $C_L$ and let $n$ be the number of connected components of $C_L$.  Then the number of valid MV-assignments of $C$ (that locally fold flat) is $2^n$.
\end{theorem}

\subsection{Example: square twist tessellations}

\begin{figure}[htp]
\centerline{\includegraphics[scale=.45]{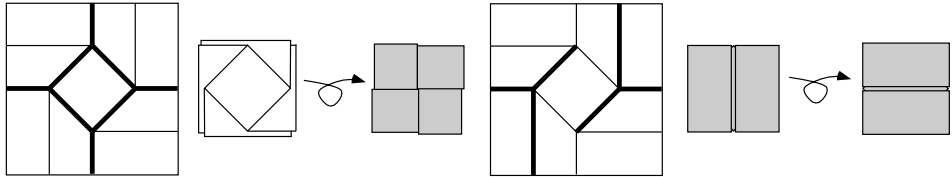}}
\caption{Two different MV-assignments of the square twist. Bold lines are mountain creases, and non-bold creases are valleys.}\label{fig:sqtwist}
\end{figure}

A {\em square twist} is an origami maneuver that literally twists a square in the paper $90^\circ$ creating four perpendicular pleats radiating from the square as it does so.  Examples are shown in Figure \ref{fig:sqtwist}, where one can see that different MV-assignments for the square twist are possible.  See \cite{H1} for a more detailed exposition on the number of ways to flat-fold a square twist.

Square twists are interesting because they easily tessellate.  That is, we can arrange square twists in a grid on a single sheet of paper, making the pleats of the twists line up and taking mirror-images of the square twist as needed.  Then the whole tessellated crease pattern should be able to fold flat (if done with care).  This is an example of an {\em origami tessellation}, a genre of origami that has become quite popular.  See \cite{Gjerde} for more examples.

\begin{figure}
\centerline{\includegraphics[scale=.7]{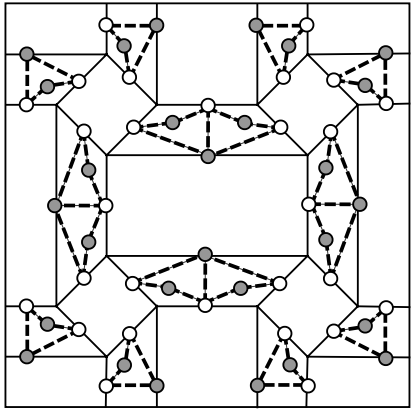}}
\caption{A $2\times 2$ tessellation of square twists (solid lines) with the origami line graph superimposed (round dots and dashed lines).}\label{fig:sqtwistcolor}
\end{figure}

We let $S(m,n)$ denote the crease pattern of an $m\times n$ square twist tessellation grid.  The $S(2,2)$ case is shown in Figure \ref{fig:sqtwistcolor}, along with its origami line graph $S(2,2)_L$.  Notice how at each vertex of this crease pattern we may employ the Big-Little-Big Lemma at the $45^\circ$ angle.  Also, by Maekawa's Theorem, the creases surrounding the $135^\circ$ angle at each vertex must have the same MV parity.  This explains the structure of $S(2,2)_L$ shown in Figure \ref{fig:sqtwistcolor}.

The general square twist origami line graph $S(m,n)_L$ will have four components at each square twist.  The $2(m+n)$ components on the border of the paper will only be touching one square twist; the others will be touching two.  This gives $(4mn-2(m+n))/2+2(m+n) = 2mn + m + n$ connected components of $S(m,n)_L$.

\begin{theorem}\label{thm:sqtwist}
The number of valid MV-assignments of $S(m,n)$ that locally fold flat is $2^{2mn+m+n}$.
\end{theorem} 

\begin{figure}
\centerline{\includegraphics[scale=.6]{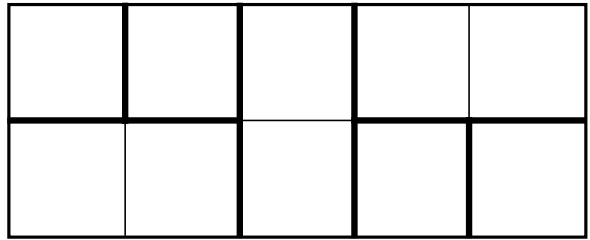}}
\caption{A $2\times 5$ stamp-folding MV assignment that is impossible to fold.}\label{fig:imposs}
\end{figure}

We conjecture that all of these MV-assignments of $S(m,n)$ are globally flat-foldable as well, but Theorem \ref{thm:sqtwist} should be viewed as an upper bound for global flat-foldability.  In fact, it is quite possible for a crease pattern to be flat-foldable at every vertex and to have no mountain-valley contradictions that the origami line graph would detect, yet still be unfoldable.  Figure \ref{fig:imposs} shows an example of this; see \cite{Justin} or \cite{GHull} for details on why it fails to globally fold flat.

\section{The Miura-ori}

\begin{figure}
\centerline{\includegraphics[scale=.4]{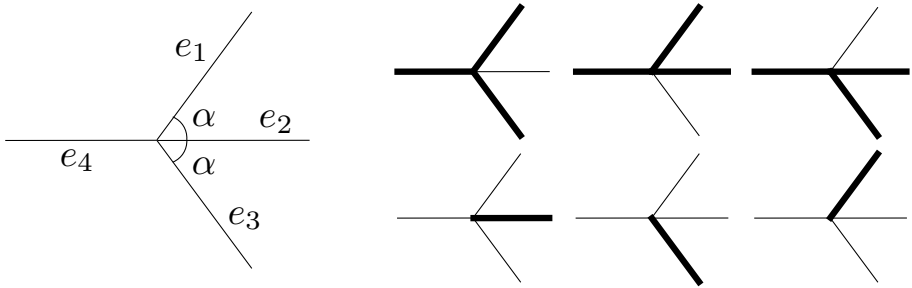}}
\caption{A flat-foldable vertex whose MV restrictions are not captured by the origami line graph of Definition \ref{def:origami line graph}.}\label{fig:Miuravertex}
\end{figure}

There are many interesting crease patterns for which Theorem \ref{thm:2color} is not applicable.  Any crease pattern whose origami line graph, as defined above, is inadequate for capturing the MV relationships between the creases will fall into this category.  A simple example can be seen in a single, degree 4 flat-foldable vertex  with two congruent acute angles adjacent, as shown in Figure \ref{fig:Miuravertex}.  Note that by Maekawa's Theorem we need three mountains and one valley (or vice-versa), and the crease labeled $e_4$ cannot be the sole valley (or the sole mountain).  Otherwise the two acute angles $\alpha$ would have to wrap around and contain the two obtuse angles $180^\circ-\alpha$, which is impossible without the paper ripping or forming new creases.  Thus, the only valid MV-assignments for such a vertex are those shown on the right side of Figure \ref{fig:Miuravertex}.  Notice that the creases labeled $e_1$, $e_2$, and $e_3$ switch in pairs from having the same to having different MV-assignments, which means the there would be no edges between these creases in the origami line graph.  Nonetheless, there are MV restrictions between $e_1$, $e_2$, and $e_3$, and thus the origami line graph will not capture these restrictions.  Other means must be used to count the number of valid MV-assignments for such vertices.

The vertex shown in Figure \ref{fig:Miuravertex} is exactly the vertex that is tessellated in the classic Miura-ori crease pattern, also known as the Miura map fold \cite{Miura1}.  The Miura-ori has attracted considerable attention over the past 30 years for its applications in engineering and nature \cite{Maha1, Wei}.  

Recently Ginepro and the author performed an analysis of MV-assignments for Miura-ori crease patterns that consist of an $m\times n$ grid of parallelograms.  This led to a bijection between the number of locally flat-foldable $m\times n$ Miura-ori MV-assignments and the number of ways to properly 3-vertex color an $m\times n$ grid graph with one vertex pre-colored.  We summarize this bijection here and refer the reader to \cite{GHull} for details of the proof.

\begin{figure}
\centerline{\includegraphics[scale=.6]{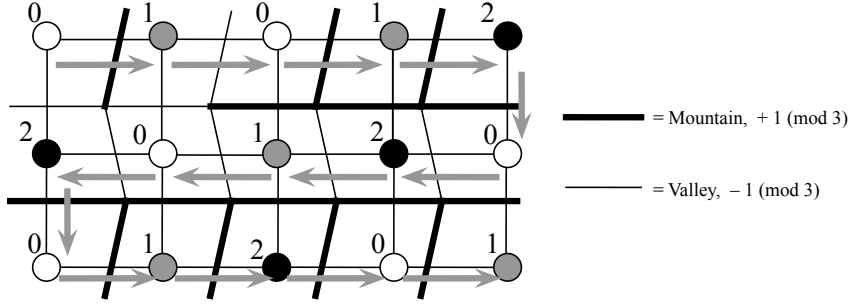}}
\caption{How we go from a locally flat-foldable Miura-ori MV-assignment to a proper 3-vertex coloring of a grid graph, and vice-versa.}\label{fig:Miurabij}
\end{figure}

The idea of the bijection is illustrated in Figure \ref{fig:Miurabij}.  Imagine overlying the $m\times n$ grid graph (with $m$ rows and $n$ columns of vertices) on top of the $m\times n$ Miura-ori so that each grid graph vertex is in the center of a parallelogram.  (In graph theory terms, the grid graph is the planar dual to the Miura-ori crease pattern, ignoring the outside face.)  The Miura-ori crease pattern should be oriented so that the top row of vertices are all ``pointing left," i.e. so that the crease $e_4$ in Figure \ref{fig:Miuravertex} is to the left of the upper-left Miura-ori vertex.  We also use as our grid graph vertex colors the integers mod three (that is, the elements of the group $\mathbb{Z}_3$) and we assume that the upper-left vertex in the grid graph gets color 0.

We then follow a zig-zag path on the $m\times n$ grid graph from the upper-left vertex, across the top row to the upper-right vertex, then down one vertex, then across the second row to the left, then down one vertex, then across to the right again, and so on.  This path is illustrated by the grey arrows in Figure \ref{fig:Miurabij}.  

We use this path to establish our bijection, which we will now describe.  Let the vertices of the $m\times n$ grid graph be denoted $v_1$, $v_2$, ..., $v_{mn}$ in the order that they are encountered on the zig-zag path.  Let $c_i$ denote the Miura-ori crease line between vertices $v_{i-1}$ and $v_i$ in the superimposed grid graph.  Our MV-assignment for the Miura-ori crease pattern $C=(V_C,E_C)$ will be $\mu:E_C\rightarrow \{-1,1\}$ and our 3-coloring of the grid graph $G=(V_G,E_G)$ will be $c:V_G\rightarrow \mathbb{Z}_3$.

{\bf To convert from $\mu$ to $c$:}  Let $c(v_1)=0$ and then recursively define
$$c(v_i) = c(v_{i-1}) +\mu(c_i)\ (\mbox{where addition is mod } 3).$$

{\bf To convert from $c$ to $\mu$:}  For creases $c_i$ between vertices $v_{i-1}$ and $v_i$ in the grid graph, define
$$\mu(c_i) = \left\{ \begin{array}{cl}
1 & \mbox{if }c(v_i)-c(v_{i-1})\equiv 1\ (\mbox{mod } 3) \\
-1 & \mbox{if }c(v_i)-c(v_{i-1})\equiv 2\ (\mbox{mod } 3).\end{array}\right.$$
For the other creases, let $d_i\in E_C$ be the crease directly below vertex $v_i$ and above $v_j$ in the superimposed grid graph.  Then define
$$\mu(d_i) = \left\{ \begin{array}{cl}
1 & \mbox{if }c(v_i)-c(v_j)\equiv 1\ (\mbox{mod } 3) \\
-1 & \mbox{if }c(v_i)-c(v_j)\equiv 2\ (\mbox{mod } 3).\end{array}\right.$$

It should not be immediately obvious to the reader that these work.  That is, converting from $\mu$ to $c$ will create a proper coloring along the zig-zag path, but a proof is required to guarantee that the coloring will be proper along $G$'s edges not in the zig-zag path.  It turns out, however, that the MV restrictions of the Miura-ori crease pattern are exactly what is needed to ensure that $c$ will be a proper coloring across these other edges.  Similarly, when constructing $\mu$ from $c$ as defined above, one needs to prove that the resulting MV assignment is locally flat-foldable.  In other words, every vertex in the Miura-ori $C$ needs to look like one of the six possibilities in Figure \ref{fig:Miuravertex} under $\mu$.  Proofs of these details are omitted here for space considerations and can be found in \cite{GHull}.  An interesting application of this bijection to a further study of the Miura-ori can be found in \cite{Hull5}.

One interesting consequence of this bijection between MV-assignments of the Miura-ori and grid graph 3-vertex colorings is that one can then use results from graph theory to gain insight into the corresponding MV-assignment counting problem.  Counting 3-colorings of grid graphs is not a completely solved problem, although the numbers generated by counting such colorings is sequence A078099 in the On-Line Encyclopedia of Integer Sequences (\url{http://oeis.org}).  Under this sequence's encyclopedia entry there is information on the transfer matrix for generating these numbers, which can thus be used to count the locally flat-foldable MV-assignments of an $m\times n$ Miura-ori.

Furthermore, in 1967 Lieb proved that this same grid graph coloring problem is equivalent to enumerating the number of states in an antiferroelectric model for two-dimensional ice lattices, otherwise known as the {\em square ice model} \cite{Lieb}.  Lieb further showed that in a grid graph with $N$ vertices for $N$ very large (say, on the order of $10^{23}$, which corresponds to the number of atoms one might have in a piece of ice), we will have 
$$(4/3)^{3N/2}$$
ways to properly 3-vertex color the grid graph with one vertex pre-colored.  Because of our bijection, this means that the number of ways to locally fold flat a Miura-ori crease pattern with $N$ parallelograms will be approximately $(4/3)^{3N/2}$ for $N$ very large.  Perhaps more importantly, this establishes a relationship between counting valid MV-assignments of origami crease patterns and Ising spin models in physics.

\section{Conclusion}

Counting the number of ways in which a crease pattern can be folded has been a challenging area of origami-mathematics.  Counting valid MV-assignments in particular has seen little progress aside from the single-vertex case \cite{Hull3}.  Developing a more general way to convert such counting problems to graph coloring would be a major breakthrough in the area.  The Ginepro-Hull bijection with Miura-ori crease patterns is very promising in this regard.

However, the bijection techniques used in the Miura-ori are very specific to that crease pattern.  It is not clear how one would generalize this to other crease patterns, especially those with vertices of degree larger than four.  There is plenty of further work to be done.

It is interesting to note how one could argue that the crease patterns that fall under the hypotheses of Theorem \ref{thm:2color} (that is, where the origami line graph tells us everything we need to know) are much more common than other crease patterns.  The idea is that in the configuration space of a flat-foldable vertex, the only components with non-zero volume are those that are {\em generic} in that the sequence of angles do not contain any consecutive equal angles or any other exploitable symmetries (see \cite{Hull4} for more details).  Such vertices would have instances of the Big-Little-Big Lemma present, and could thus have all of their MV restrictions captured by the origami line graph.  It stands to reason that if we were to take a ``generic" flat-foldable crease pattern (all of whose vertices are generic) then its MV restrictions might be completely described by the origami line graph.  In other words, if we were to use these ideas to define a {\em random flat-foldable crease pattern} then it could be true that almost every flat-foldable crease pattern could have its MV-assignments enumerated by the origami line graph and Theorem \ref{thm:2color}.  

Nonetheless, it is also true that most of the crease patterns that we find interesting, like the Miura-ori, possess some amount of symmetry and thus would not be generic.  Such crease pattern MV-assignments will remain more challenging to enumerate.

\section{Acknowledgements}

The author would like to thank Crystal Wang for useful discussions leading up to the preparation of this paper.  This research was supported by the National Science Foundation grant EFRI-ODISSEI-1240441 ``Mechanical Meta-Materials from Self-Folding Polymer Sheets."

\bibliographystyle{akpbib}
\bibliography{colorcounting.bib}

\end{document}